%
% This is a plain TeX file
%
 \magnification \magstep1
\font\bold=cmbx10 at 14pt

\centerline{\bold Three-Dimensional Manifolds, Skew-Gorenstein Rings }
\medskip 
\centerline{\bold  and their Cohomology.}
\bigskip
\centerline{ Jan-Erik Roos}
\centerline{Department of Mathematics}
\centerline{Stockholm University}
\centerline{SE--106 91 Stockholm, SWEDEN}
\centerline{ e-mail: {\tt jeroos@math.su.se}}
\bigskip
% \rightline{ \it Dedicated to Ralf Fr{\"o}berg and Clas L{\"o}fwall at their $65^{th}$ birthdays.}
\bigskip

\centerline{January 13, 2010}

\bigskip

\def\mysec#1{\bigskip\centerline{\bf #1}\nobreak\par}

\def\cite#1{~[{\bf #1}]}
{\openup 1\jot

\mysec{ Abstract.}
\bigskip
Graded skew-commutative rings occur often in practice. Here are two examples:
1) The cohomology ring of a compact three-dimensional manifold.
2) The cohomology ring of the complement of a hyperplane arrangement (the Orlik-Solomon algebra).
We present some applications of the homological theory of these graded skew-commutative rings. In particular we find compact oriented 3-manifolds 
without boundary for which the Hilbert series of 
the Yoneda Ext-algebra of the cohomology ring of the fundamental group
is an explicit transcendental function.
This is only possible for large first Betti numbers of the 3-manifold (bigger than -- or maybe equal to -- 11).
We give also examples of 3-manifolds where the Ext-algebra of the cohomology ring of the fundamental group
is not finitely generated.


\mysec{0. Introduction}
\bigskip
Let $X$ be an oriented compact $3$-dimensional manifold without boundary.
The cohomology ring $H=H^*(X,{\bf Q})$
is a graded skew-commutative ring whose augmentation ideal $\bar H$ satisfies ${\bar H}^4=0$.
The triple (cup) product $x \cup y \cup z =\mu(x,y,z).e  $ ,
where e is the orienation generator of $H^3$ defines a skew-symmetric trilinear form
on $H^1$ with values in $\bf Q$  and conversely, according to a theorem of Sullivan [Sul] any such form comes in this way from a 3-manifold $X$ (not unique)
whose cohomology algebra can be reconstructed from $\mu$ since by Poincar\'e duality $H^2 \simeq (H^1)^*$.
In the more precise case when $H^*$ is a also a Poincar\'e duality algebra, i.e.
the cup product $H^1 \times H^2 \rightarrow H^3$ is non-degenerate, it
follows that $H^*$ is a Gorenstein ring (cf. section 1 below).
 Such Gorenstein rings will be studied here. Any 3-manifold $M$ can be decomposed in a unique way
as a connected sum of prime 3-manifolds:
$$
M = P_1 \sharp P_2 \cdots \sharp P_k
$$
of prime manifolds $P_i$ (cf e.g. Milnor [Mil], Theorem 1), and with the exception of $S^3$ and $S^2 \times S^1$
any prime manifold is also irreducible ([Mil] and for them $\pi_2(M)=0$ [Mil, Theorem 2]. If furthermore
$\pi_1(M)$ is infinite then also the higher homotopy groups are 0 ($M$ is said to be aspherical)
so that $M$ is the Eilenberg-Maclane space
$K(\pi_1(M),1)$ and the cohomology ring of $M$ is isomorphic to the cohomology ring of the group 
$\pi_1(M)$  (this is of course also true for any 3-manifold $X$ which is aspherical.
(For most of the applications below we suppose that $H=H^*(X)$ is a Poincar\'e
duality algebra and we suppose that the base field $k$ is of characeristic 0, the preference being {\bf Q}.)
The ring $H$ has interesting homological properties
which have not yet been fully studied, and we wish to continue such a study here.
For small values of the first Betti number $b_1(X)= dim_{\bf Q}H^1(X,{\bf Q})$
 the ring $H$ is a Koszul algebra (cf. section 2 below) so that in particular the generating series
$$
P_H(z) = \sum_{i\geq 0}|Tor_i^H(k,k)|z^i = H(-z)^{-1}  \leqno(0.1)
$$
where $H(z) = 1 + |H^1(X,Q)|z + |H^2(X,Q)|z^2+z^3$,
and where $|V|$ denotes the dimension of the vector space $V$. 
 But for bigger Betti numbers many new
phenomena occur. In particular we will see that for $b_1(X)= 12$
(and maybe even for $b_1(X) =11$) there are a few  examples where
$P_H(z)$ is an explicit transcendental function (thus we are far away
from the Koszul case of formula (1)!). But maybe 11 is the best possible number here.
On the other hand for bigger $b_1(X)$ the possible series $P_H(z)$ are
rationally related to the family of series which occurred in
connection with the Kaplansky-Serre questions a long time ago [An-Gu].
But even for smaller $b_1(X)$ other strange homological properties of $H$  
occur: we will give examples (probably best possible) where $b_1(X)=11$ and the Yoneda $Ext$-algebra
$Ext^*_H(k,k)$ is not finitely generated.  
The implications of all this for $3$-manifold groups have not been
fully explored. Cf. also [Sik]. Conversely, the connection with $3$-manifolds
makes it possible to go backwards and to deduce results in the homology theory of skew-commutative algebras (and there are also related
results in the {\it commutative} case).
 Note that we are studying everything over $\bf{Q}$.
There seem to be some relations with Benson [Be] but he works
over finite fields.
Finally, let us remark that even for the special 3-dimensional case when $X$
is the boundary manifold of a line arrangement in $P^2({\bf C})$ 
we can have e.g. the same strange transcendental phenomenon as above but the prize to pay for this
is to accept even bigger $b_1(X)$. The three-manifolds that occur here are called
graphic manifolds [Co-Su] and they are aspherical if the line arrangement is not a pencil 
of lines. 
\bigskip
\mysec{1. Graded skew-Gorenstein rings and classification of trivectors}
\bigskip
Let us first recall that a local commutative Gorenstein ring was defined in [Bass] as
a ring $R$ which has a finite injective resolution as a module over itself.
In particular if $R$ is artinian this means that $R$ is injective as a module over itself.
It is also equivalent to saying that the socle of $R$ , i.e. $Hom_R(R/m,R)$ ($m$ is the maximal ideal of $R$)
is 1-dimensional over $k=R/m$. 
Things are more complicated in the noncommutative case [Fo-Gr-Rei], but if
$R$ is skew-commutative artinian the Gorenstein property is equivalent to $R$ being injective as a module
over itself (left or right --- these two conditions are equivalent --- and they are also equivalent to
saying that the left --- or right ---  socle of $R$ is one-dimensional).   
In the special case when $X$ be an oriented compact three-dimensional manifold without boundary, and
where  $H = H^*(X,{\bf Q})$
and let $R=H^*(X,{\bf Q})$ be the cohomology ring of $X$. From the preceding definition it follows that
$R$ is Gorenstein if and only if $R$ is generated by $H^1$ and $H$ is a Poincar\'e duality algebra
(we assume that $|H^1| > 1$).

We now turn to the classification of such $R$:s when $|H^1(X,{\bf Q})|\leq 8$. We will use the
classification of trivectors in $H^1$ described in paragraph 35 of the book [Gur].

For the rank 5 there is only one trivector (denoted by III in [Gur], page 391) and which can be denoted
by (here $e_1,e_2,e_3,e_4,e_5$ is a basis for $H^1$ and the $e^i$:s are the dual basis elements in $(H^1)^*$):
$$
e^1\wedge e^4 \wedge e^3 + e^2 \wedge e^5 \wedge e^3
$$
and the corresponding Gorenstein ring
is a Koszul algebra (there is a small misprint in [Sik]) 
It is only when the ranks are $\geq 6$ that non-Koszul Gorenstein rings occur.
Let us give the details in the first nontrivial case of [Gur], namely case IV of rank 6, page 391:
$$
f=e^1 \wedge e^2 \wedge e^3 +e^3 \wedge e^4 \wedge e^5 + e^2 \wedge e^5 \wedge e^6  
$$
We want to determine all elements $g=\sum_{i<j}\mu^s_{i,j}e_i\wedge e_j$ in the exterior
algebra 
$$
E(e_1,e_2,e_3,e_4,e_5,e_6)
$$
such that
$\sum_{i<j}\mu^s_{i,j}f(e_i,e_j,e_s) = 0$
for $s=1,\ldots, 6$.
This gives, using the explicit form of $f$ and the condition that $f$ is skew-symmetric the conditions:
$$ 
\mu^3_{1,2}=1,\quad \mu^2_{1,3}=-1,\quad \mu^1_{2,3}=1 \quad,\mu^6_{2,5}=1 \quad,\mu^5_{2,6}=-1
$$
$$ 
\mu^5_{3,4}=1,\quad \mu^4_{3,5}=-1,\quad \mu^3_{4,5}=1,\quad \mu^2_{5,6}=1,
$$
and all other $\mu^s_{i,j}=0$, leading to the following ring with quadratic relations (we do not write wedge for
multiplication):
$$
R_{IV} = {E(e_1,e_2,e_3,e_4,e_5,e_6) \over (e_1e_2-e_4e_5,e_1e_3+e_5e_6,e_2e_6+e_3e_4,e_1e_4,e_1e_5,e_1e_6,e_2e_4,e_3e_6,e_4e_6)} 
$$
But this quotient ring $R_{IV}$ has Hilbert series $R(z)=1+6z+6z^2+2z^3$. This means that we have to find a last
cubic relation. One finds that the cube of the maximal ideal of $R$
is generated by $e_2e_5e_6$ and $e_2e_3e_5$.
But $f(e_2,e_3,e_5)=0$ and $f(e_2,e_5,e_6)=1$ so that the corresponding 
Gorenstein ring is $G_{IV}=R/e_2e_3e_5$ and we will see that this Gorenstein ring is {\it not} a Koszul algebra.
For the case V of {\it loc.cit.} the $f$ is given by
$$
f = e^1 \wedge e^2 \wedge e^3 + e^4 \wedge e^5 \wedge e^6
$$
leading in the same way to the ring with quadratic relations:
$$
R_V ={E(e_1,e_2,e_3,e_4,e_5,e_6) \over(e_1e_4,e_1e_5,e_1e_6,e_2e_4,e_2e_5,e_2e_6,e_3e_4,e_3e_5,e_3e_6)}
$$
which also has Hilbert series $R(z)=1+6z+6z^2+2z^3$, but in this case we have to divide by
$e_1e_2e_3-e_4e_5e_6$ to get the Gorenstein ring $G_V=R_V/(e_1e_2e_3-e_4e_5e_6)$
which is, as we will see below, not a Koszul algebra either.
\bigskip
In [Gur, pages 393-395] the classification of 3-forms of rank 7 is given as the 5 cases VI,VII,VIII,IX,X and those forms
of rank 8  are described as the 13 cases XI,XII,...,XXIII.
The classification of 3-forms of rank 9 are given in [Vin-El].
We will describe the homological behaviour of the corresponding Gorenstein rings in the section 4.
 \mysec{2. Calculating the Koszul dual of the Gorenstein ring associated to a 3-form.}
\bigskip
Let $R$ be any finitely presented ring (connected k-algebra) generated in degree 1 and having quadratic relations.
It can be described as the quotient $T(V)/(W)$, where $T(V)$ is the tensor algebra on a (finite-dimensional)
$k$-vector space $V$ , placed in degree 1 and $(W)$ is the ideal in $T(V)$, defined 
by a sub-vector space $W \subset V\otimes_k V$. The Yoneda Ext-algebra of R is defined by
$$
Ext^*_R(k,k) = \bigoplus_{i\geq 0}Ext^i_R(k,k) 
$$ 
where $k$ is an $R$-module in the natural way and where the multiplication is the Yoneda product.
The sub-algebra of $Ext^*_R(k,k)$ generated by $Ext^1_R(k,k)$ is called the Koszul dual of $R$ and it is
denoted by $R^!$. It can be calculated as follows: consider the inclusion map $W \rightarrow  V\otimes_k V$.
Taking k-vector space duals (denoted by $W^*$ and $(V\otimes_k V)^*$ we get the exact sequence
($W^{\perp}$ = those linear $f: V\otimes_k V \rightarrow k$ that are $0$ on $W$):
$$
   0 \leftarrow W^* \leftarrow (V\otimes_k V)^* \leftarrow W^{\perp} \leftarrow 0
$$ 
Now $R^! = T(V^*)/(W^{\perp})$  (we have used that  $(V\otimes_k V)^* =V^*\otimes_k V^*$).

For all this cf. [L\"of]. Note that $(R^!)^!$ is isomorphic to $R$.
Note also that if $R$ also has cubic relations and/or higher relations then
$Ext^*_R(k,k)$ is still defined, and the subalgebra generated by $Ext^1_R(k,k)$ is still given by the formula
 $T(V^*)/(W^{\perp})$ where $W$ is now only the ``quadratic part'' of the relations of $R$. In particular
$(R^!)^!$ is now only isomorphic to $T(V)$ divided by the quadratic part of the relations.
Not also that in general, if $R$ is skewcommutative then $Ext^*_R(k,k)$ is a cocommutative Hopf algebra
which is the enveloping algebra of a graded Lie algebra, and $R^!$ is a sub Hopf algebra which is the
enveloping algebra of a smaller graded Lie algebra. (All this is also true if $R$ is commutative,
but now the Lie algebras are super Lie algebras.)
In [L\"of] Corollary 1.3, pages 301-302 there is a recipe about how to calculate the Koszul dual of an
algebra with quadratic relations.
Applying this to the case of $R_{IV}$ of the previous section we find that $R_{IV}^!$ is the algebra
$$
{k<X_1,X_2,X_3,X_4,X_5,X_6> \over
([X_1,X_2]+[X_4,X_5],[X_3,X_4]-[X_2,X_6],[X_1,X_3]-[X_5,X_6],[X_2,X_3],[X_2,X_5],[X_3,X_5])}
$$
where $[X_i,X_j]=X_iX_j-X_jX_i$ is for $i<j$ the Lie commutator of $X_i$ and
$X_j$.
In general if one starts with a (skew)-commutative algebra $A$ it is often easy to calculate the Hilbert series
of $A$. But calculating the Hilbert series of $A^!$ is often very difficult and this series can even be a
transcendental function.
But in this case it is rather easy: we get that $R_{IV}^!(z) = 1/(1-6z+6z^2-2z^3)$ so that 
$R_{IV}(-z)R_{IV}^!(z) = 1$ and we even have (the in general strictly stronger assertion ([Ro 4],[Po])
that $R_{IV}$ is a Koszul algebra.
(For the definitions and equivalent characterizations of Koszul algebras we refer to
[L\"of, p. 305, Theorem 1.2].)  A similar result for $R_V$ holds true, and in this case we
can directly apply a result of Fr\"oberg [Fr\"o], since $R_V$ has quadratic monomial relations. 
But neither the Gorenstein quotients $G_{IV}$ nor $G_V$ are Koszul algebras and in the next section
we will see how to relate their Hilbert series and the corresponding Hilbert series of their Koszul duals
to their two-variable Poincar\'e-Betti series 
$$
P_G(x,y)=\sum_{i,j} |Tor_{i,j}^G(k,k)|x^iy^j
$$
But we will first deduce a few results about explicitly calculating the Gorenstein ring and its 
Koszul dual associated to a given 3-form.

An old result of Macaulay gives a nice correspondence between commutative artinian graded Gorenstein algebras having
socle of degree $j$
of the form $k[x_1,x_2,\ldots,x_n]/I$ and homogeneous forms of degree $j$ in the dual of $k[x_1,x_2,\ldots,x_n]$
({\it cf. eg.} Lemma 2.4 of [ElK-Sri]). Here is a skew-commutative version 
(here announced only for socle  degree 3) which we have not been able to
find in the literature (I thank Antony Iarrobino who suggested that such a result should be true):

PROPOSITION 2.1 (``Skew''-Macaulay).- Let $E=E[X_1,\ldots, X_n]$ be the exterior algebra in $n$ variables of degree 1
over a field of characteristic $0$. Let $E[\partial/\partial X_1,\ldots,\partial/\partial X_n]$
be the exterior algebra of skew-derivations of $E$. To each homogeneous element $F$ of $E^*$ of degree 3 we associate
the elements $h$ in $E$, such that $h(\partial/\partial X_1,\ldots,\partial/\partial X_n)(F) = 0$
Let $I$ be the ideal of these elements in $E$ and $R=E/I$. Then $R$ is a skew-Gorenstein ring, where $m^4 = 0$.
Conversely, given such a Gorenstein ring $R=E/I$, we can recover $F$ as a generator of the degree 3 part of
the ideal in $E^*$
consisting of those $K$  such that  $K(\partial/\partial X_1,\ldots,\partial/\partial X_n)f =0$
for all f in $I_3$.
\bigskip
Note that this is essentially the reasoning we used when we calculated the degree 2 part of $I$ in
the case IV above. Combining this with the L\"ofwall description of the calculation of $R^!$ given above
we finally arrive at the following useful Theorem-Recipe to calculate the Koszul dual $G^!$ of the Gorenstein
ring associated to a skew 3-form:
\bigskip   
THEOREM-RECIPE 2.1 --- Let $\Psi(e^1,e^2,\ldots,e^n)$ be a skew-symmetric 3-form of rank $n$, $X$ one of the 3-manifolds
giving rize to $\Psi$ and $G = H^*(X,{\bf Q})$. The Koszul dual $G^!$ of $G$ is obtained as follows:
Calculate the $n$ skew-derivaties of $\Psi$ with respect to the variables $e^1,e^2,\ldots,e^n$
Then
$$
G^! \simeq {k<X_1,X_2,\ldots,X_n>\over (q_1,q_2,\ldots,q_n)}
$$
where $k<X_1,X_2,\ldots,X_n>$ is the free associative algebra generated by the variables $X_i$
that are {\it dual} to the $e^i$:s in $\Psi$ and where the $q_i$:s are obtained by replacing
each quadratic element $e^s\wedge e^t$ ($s < t$ in the $\partial \Psi /\partial e^i$ by the commutator
$[X_s,X_t] =X_sX_t-X_tX_s$ in $k<X_1,X_2,\ldots,X_n>$, for $i = 1,\ldots n$.

Note that if $G^!$ is given it is easy to calculate ``backwards'' the ring with quadratic relations
$(G^!)^!$. After that it is easy to find the extra cubic relations we should divide with to get $G$.
Using this theorem-recipe we finally find:
\bigskip
THEOREM 2.2 --- a) The double Koszul duals $(G^!)^!$ of the Gorenstein rings corresponding trivectors
of rank $= 7$ are Koszul algebras in the cases VI,VII, VIII,IX  and X and therefore $G^!$ is a Koszul
algebra in cases VI, VII, VIII, IX, X. (But they are not isomorphic). The corresponding Gorenstein algebras
$G_{VI},G_{VII},G_{VIII},G_{IX},G_{X}$ are also Koszul algebras.
\bigskip
b) When it comes to trivectors of rank 8, i.e. the cases
XI,XII,$\ldots$ XXIII the situation is more complicated already from the
homological point of view:
Indeed in case XI the double Koszul dual $(G^!)^!$ is already a 
Gorenstein ring with quadratic relations, thus equal to $G$ and $G^!$ is
{\it not} a Koszul algebra.
But the case XII treated explicitly below is slightly different and not a Koszul algebra either. 
The cases XIII, XIV and XV are Koszul algebras.
But case XVI is as above: $(G^!)^!$ is a Koszul algebra, but we have
to  divide out a cubic form to get $G_{XVI}$ which is therefore not Koszul.
Finally the algebras corresponding to the cases XVII,XVIII,XIX,XX,XXI,XXII and XXIII are Koszul algebras (but not isomorphic). 
\bigskip 
EXAMPLE 2.1 Here is a use of the THEOREM-RECIPE 2.1:
Let us consider the case XII of rank 8.
Here the 3-form $f_{XII}=\Psi$ is
$$
\Psi(e^1,e^2,\ldots,e^8)=e^5\wedge e^6 \wedge e^7 +e^1 \wedge e^5 \wedge e^4 + e^2\wedge e^6 \wedge e^4 
+e^3\wedge e^7 \wedge e^4+ e^3\wedge e^6 \wedge e^8 
$$ 
We calculate the 8 partial derivaties of this skew form (we do not write out the $\wedge$ sign):
$$
{\scriptstyle \partial \Psi /\partial e^1 = e^5e^4,\quad \partial \Psi /\partial e^2 = e^6e^4, 
\quad \partial \Psi /\partial e^3 = e^7e^4+e^6e^8,\quad \partial \Psi /\partial e^4 = e^1e^5+e^2e^6+e^3e^7,}
$$
$$
{\scriptstyle \partial \Psi /\partial e^5 = e^6e^7-e^1e^4,\quad \partial \Psi /\partial e^6 = -e^5e^7-e^2e^4-e^3e^8,\quad
\partial \Psi /\partial e^7 = e^5e^6-e^3e^4, \quad \partial \Psi /\partial e^8 = e^3e^6} 
$$
leading to $G^! = k<X_1,X_2,X_3,X_4,X_5,X_6,X_7,X_8>$ divided by the ideal
$$
 \scriptstyle {( [X_5,X_4],[X_6,X_4],[X_7,X_4]+[X_6,X_8],
[X_1,X_5]+[X_2,X_6]+[X_3,X_7],[X_6,X_7]-[X_1,X_4],}
$$
$$
 \scriptstyle {
[X_5,X_7]+[X_2,X_4]+[X_3,X_8],[X_5,X_6]-[X_3,X_4],[X_3,X_6])}
$$
and the Hilbert series $1/G^!(z)= 1-8z+8z^2-z^3-z^4$. Furthermore one sees that $(G^!)^!$ has Hilbert series
 $1+8z+8z^2+z^3$  so that we do not have to divide by cubic elements to get $G_{XII}$ which is non-Koszul.
As a matter of fact we have a generalization of a formula of L\"ofwall (cf. next section), giving that
$$
{1\over P_{G_{XII}}(x,y)} = (1+1/x)/G_{XII}^!(xy)-G_{XII}(-xy)/x
$$
so that 
$$
{1\over P_{G_{XII}}(x,y)} = 1-8xy+8x^2y^2-x^3y^3-x^3y^4-x^4y^4
$$
leading to {\it e.g.}  $Tor_{3,4}^{G_{XII}}(k,k)$ being 1-dimensional.
The other cases I-XXIII are treated in a similar manner: sometimes $G^!(z)$
can be directly calculated since we have a finite Groebner basis for the
non-commutative ideal and in all cases using the Backelin et al programme BERGMAN [Ba].
\mysec{3. A generalization of a formula of L\"ofwall to Gorenstein rings with $m^4=0$}
\bigskip
In his thesis [L\"of] Clas L\"ofwall proved in particular that if $A$ is any graded connected
algebra with ${\bar A}^3 =0$ then the double Poincar\'e-Betti series is given by the formula
$$
{1\over P_A(x,y)} = (1+1/x)/A^!(xy)-A(-xy)/x , \leqno(3.1)
$$
In [Ro 1] we have an easy proof of this in the case when $A$ is commutative (works also 
in the skew-commutative case [Ro 2]) using the fact that in these cases $A^!$ is a sub-Hopf algebra of the 
big Hopf algebra $Ext^*_A(k,k)$ and according to a theorem of Milnor-Moore this big Hopf algebra is
free as a module over $A^!$. We now show
\bigskip
THEOREM 3.1.--- The formula (3.1) is true when $(R,m)$ is a graded (skew-) commutative Gorenstein ring
with $m^4=0$.
\medskip
PROOF: First we note that Avramov-Levin have proved [Av-Le] in the commutative case that the natural map
$$
 R \rightarrow R/soc(R)
$$
is a Golod map (cf. section 5 for the skew-commutative case that we use here).
 Now the socle of $R$ is $m^3$ which is $s^{-3}k$
Therefore we have
$$
P_R/soc(R)(x,y) = {P_R(x,y)\over 1-xy(P_R^{R/soc(R)}(x,y)-1)} = {P_R(x,y)\over 1-x^2y^3P_R(x,y)} \leqno(3.2)
$$
But $R/m^3$ is a local ring where the cube of the maximal ideal is 0.
Thus the formula of L\"ofwall can be applied, and we obtain using that $(R/m^3)^! = R^!$ the formula
$$
{1\over P_R/m^3(x,y)} = (1+1/x)/R^!(xy)-(R/m^3)(-xy)/x
$$
which combined with (3.2) gives (3.1) since $R(z) = R/m^3(z) +z^3$.
\bigskip
REMARK 3.1.--- Thus it is clear that we can not obtain cases where the ring
$H^*(X,{\bf Q})$ has ``bad'' homological properties if $b_1(X) \leq 8$
We therefore study the case when $b_1(X)=9$. In this case there is a classification of the tri-vectors
by E.B. Vinberg and A.G. Elashvili [Vin-El].
But even in this case it seems impossible to get ``exotic'' $H^*(X,{\bf Q})$.
In fact we have e.g. by the procedure above analyzed all cases
in the Table 6 (pages 69-72)in [Vin-El] (those cases that have a * in loc.cit. correspond to $b_1(X)<9$ and
they have already been treated). 
\bigskip
THEOREM 3.2.--- For the 3-forms of rank 9 of Table 6 (pages 69-72) of [Vin-El] we have that the corresponding
$1/R^!(z)=1-9z+9*z^2-z^3$, and the corresponding Gorenstein ring is a Koszul algebra in all cases except
\medskip
a) the cases 79, 81, 85 where $1/R^!(z)= 1-9z+9z^3-3z^3$
\medskip
b) the case 83    where $1/R^!(z)= 1-9z+9z^3-z^3-5z^4+4z^5-z^6$
\bigskip
PROOF.- The proof is by using the THEOREM-RECIPE 2.1  above.
We only indicate what happens in the exceptional cases a) and b).
In case 79 in a) the 3-form is (we do not out the $wedge$ for multiplication):
$$
f_{79}=e^1e^2e^9+e^1e^3e^8+e^2e^3e^7+e^4e^5e^6
$$
We take the 9 partial derivatives
$\partial f_{79} /\partial e^i$ for $i=1\ldots 9$
and we obtain:
$$
e^2e^9+e^3e^8,\quad -e^1e^9+e^3e^7,\quad -e^1e^8-e^2e^7,\quad e^5e^6,\quad -e^4e^6,\quad e^4e^5,
\quad e^2e^3,\quad e^1e^3,\quad e^1e^2
$$
leading to $R^! = $
$$
 {{\scriptstyle k<X_1,X_2,X_3,X_4,X_5,X_6,X_7,X_8,X_9>} \over {\scriptscriptstyle ( [X_2,X_9]+[X_3,X_8],
-[X_1,X_9]+[X_3,X_7],[X_1,X_8]+[X_2,X_7],[X_5,X_6],[X_4,X_6],[X_4,X_5],[X_2,X_3],[X_1,X_3],[X_2,X_3])}}
$$
Now it turns out that the Gr\"obner basis of the ideal above is finite
and in degree $\leq 2$. This proves that $R^!$ is Koszul and has
$R^!(z) = 1-9z+9z^3-3z^3$, since the commutative ring $(R^!)^!$ has
Hilbert series $1+9z+9z^2+3z^3$.
The ideal $m^3$ is generated by the three elements ($e_1e_2e_9,e_4e_5e_6,e_1e_2e_3$)
and
$$
f_{79}(e_1,e_2,e_9)=1,\quad f_{79}(e_4,e_5,e_6)=1 \quad {\rm and}\quad f_{79}(e_1,e_2,e_3)=0
$$ 
It follows that if we divide out by the two extra elements $e_1e_2e_9-e_4e_5e_6,e_1e_2e_3$ we do have the non-Koszul
Gorenstein ring $G_{79}$ we are looking for. The rings $G_{81}$ and $G_{85}$ in a) are treated in the same way.
They correspond to:
$$
f_{81}= e^1e^2e^9+e^1e^3e^8+e^1e^4e^6+e^2e^3e^7+e^2e^4e^5+e^3e^5e^6
$$
and
$$
f_{85}=e^1e^2e^9+e^1e^3e^5+e^1e^4e^6+e^2e^3e^7+e^2e^4e^8
$$ 
\medskip
In case 83 of b) the form 
$$
f_{83}=e^1e^2e^9+e^1e^3e^5+e^1e^4e^6+e^2e^3e^7+e^2e^4e^8+e^3e^4e^9 
$$
and its 9  skew-derivaties lead to the quotient:
$$
\scriptstyle {k<X_1,X_2,X_3,X_4,X_5,X_6,X_7,X_8,X_9>}
$$
divided by the ideal
$$
\scriptstyle { ( [X_2,X_9]+[X_3,X_5]+[X_4,X_6],
-[X_1,X_9]+[X_3,X_7]+[X_4,X_8],-[X_1,X_5]-[X_2,X_7]+[X_4,X_9],}
$$
$$\scriptstyle {[X_1,X_6]+[X_2,X_8]+[X_3,X_9],
[X_1,X_3],[X_1,X_4],[X_2,X_3],[X_2,X_4],[X_1,X_2]+[X_3,X_4])}
$$
In this case the Hilbert series of $(R^!)^!$ is
$1/(1-9z+9z^2-z^3-5z^4+4z^5-z^6)$,
and the ring $(R^!)^!$ has Hilbert series:
$1+9z+9z^2+z^3$ so that $G_{83} = (R^!)^!$ which is a Gorenstein non-Koszul algebra and the THEOREM 3.2 is proved.
\bigskip
 
Thus to obtain examples of $X$ where the homological properties of the cohomology algebra $H^*(X,Q)$ are complicated
we need to study the cases where the dimension of $H^1(X,{\bf Q})$ is $\geq 10$.
But the trivectors of rank 10 have not been classified, and even the cases
of Sikora [Sik],Theorem 1 which uses for any simple Lie algebra $g$ 
the trilinear skew-symmetric form $\Psi_g(x,y,z) =\kappa(x,[y,z])$ where $\kappa$ is the Killing form of $g$
do not seem to given anything exotic from our point of view, at least when the dimension of the Lie algebra is 10.
 We therefore need a new way of constructing 
Gorenstein rings corresponding to trivectors of rank $\geq 10$, but defined
in another way.
This will be done in the next section.
\mysec{4. New skew-Gorenstein rings.}
Let $(R,m)$ be any ring with $m^3 = 0$ which is the quotient of the
exterior algebra $E(x_1,\ldots,x_n)$ by homogeneous forms of degrees $\geq 2$
and let $I(k)$ be  the injective envelope of the residue field $k=R/m$ of $R$.
Then $G=R \propto I(k)$ is a skew-Gorenstein ring with maximal ideal
$n = m \oplus I(k)$ satisfying $n^4=0$, which therefore corresponds to a 
trivector of rank $|m/m^2|+|I(k)/mI(k)|=|m/m^2|+|m^2|$ and
according to the Sullivan theory it comes from a 3-manifold $X$ with
the dimension of $H^1(X,{\bf Q})$ being equal to $|m/m^2|+|m^2|$,
 It turns out that the homological properties of $G$ are closely
related to those of the smaller ring $R$.
Indeed we have the following general result:
\bigskip
THEOREM 4.1 (Gulliksen).--- Let $R$ be (skew-)commutative and $R \propto M$ be the trivial extension of $R$ with the $R$-module $M$.
Then we have an exact sequence of Hopf algebras:
$$
k \rightarrow T(s^{-1}Ext^*_R(M,k)) \rightarrow Ext^*_{R \propto M}(k,k) \rightarrow Ext^*_R(k,k) \rightarrow k
$$
where $T(s^{-1}Ext^*_R(M,k))$ is the free algebra on the graded vector space $s^{-1}Ext^*_R(M,k)$ and where
the arrow $Ext^*_{R \propto M}(k,k) \rightarrow Ext^*_R(k,k)$ is a split epi-morphism (there is a splitting
ring map $R\propto M \rightarrow R$.
In particular we have the formula for Poincar\'e-Betti series
$$
P_{R\propto M}(x,y) = {P_R(x,y)\over 1-xyP_R^M(x,y)}
$$
\bigskip
Now turn to the artinian case and
 assume that $M$ is finitely generated and  let $\tilde M = Hom_R(M,I(k)$
be the Matlis dual of $M$ [Mat] (if $M$ as an $R$-module is a vector space over $k$, then this is the  
the ordinary vector space dual)
In this case we have a well-known formula (cf e.g. Lescot [Les 2], Lemme 1.1) 
$$
 Ext^*_R(M,k) \simeq  Ext^*_R(k,\tilde M)
$$
In particular if $M = I(k)$ then $\tilde M = Hom_R(I(k),I(k)) = R$ [Mat]
so that 
$$
 Ext^*_R(I(k),k)\simeq  Ext^*_R(k,R)
$$
Thus, if we use the Theorem on $R \propto I(k)$ we are led to the study of the Bass series of $R$,
i.e. the generating series in one or two variables of $Ext^*_R(k,R)$.
It turns out that in many (most) of the cases we study here , the Bass series of $R$ divided by
the Poincar\'e-Betti series $P_R$ is a very nice explicit polynomial (for more details
about this -- the B{\o}gvad formula - we refer to the last section of this paper.
Thus if we want strange homological properties of the Ext-algebra of the Gorenstein ring
$R \propto I(k)$ i.e the Ext-algebra of the cohomology ring of the corresponding fundamental group
$\pi_1(X)$ of the 3-manifold $X$ corresponding to the Gorenstein ring $R \propto I(k)$
we only have to find $(R,m)$ with $m^3=0$ with strange properties.
\bigskip
COROLLARY 4.1 ---  
In the case when $R$ is the cohomology ring (over ${\bf Q}$ of the complement
of a line arrangement $L$ in $P^2({\bf C})$ i.e. the Orlik-Solomon
algebra of $L$, the 3-manifold $X$ correponding to the
Gorenstein ring $R\propto I(R/m)$ can be chosen as the boundary manifold
of a tubular neighborhood of $L$ in $P^2({\bf C})$ [Co-Su]
In this case we found in [Ro 2] an arrangement where the Ext-algebra of $R$
was not finitely generated. Then the Ext-algebra of the cohomology ring of X that corresponds to
$R \propto I(k)$ can not be finitely generated either since by THEOREM 4.1  it is mapped onto $Ext^*_R(k,k)$
In this case the dimension of $H^1(X,{\bf Q})$ is 12.
Also in [Ro 2] we also found two cases of arrangements: the maclane arrangement
and the macleas arrangement where the corresponding $R^!(z)$ is a transcendental function, and this leads
to two cases where $H^1(X,{\bf Q})$ is of dimension 20, resp. 21. 
In order to press down this dimension to 12 and maybe to 11, I have to use
some of my earlier results.
In our paper [Ro 3] describing the homological properties of quotients
of exterior algebras in 5 variables by quadratic forms (there are 49 cases found), we have
found 3 cases (cases 12, 15 and 20) where $R^!(z)$ is proved to be transcendental and is
explicitly given, and 3 other cases (cases 21,22 and 33) where we conjecture that $R^!(z)$
is transcendental (but no explicit formula can be given, even in the case
33 where we have now calculated the series $R^!(z)$ up to degree 32 using Backelin's et al
programme BERGMAN (some details are given in [Ro 2], where the ``educated guess'' now has
to be abandoned):

Here is the Case 20:
$$
R_{20}={E(x_1,x_2,x_3,x_4,x_5)\over(x_1x_4+x_2x_3,x_1x_5+x_2x_4,x_2x_5+x_3x_4)}
$$
Here the Hilbert series is $R_{20}(z)=1+5z+7z^2$,
and the Koszul dual $R_{20}^!$ is given (according to the recipe we have described above
$$
k<X_1,X_2,X_3,X_4,X_5> \over {\scriptstyle([X_1,X_2],[X_1,X_3],[X_3,X_5],[X_4,X_5],[X_1,X_4]-[X_2,X_3],[X_1,X_5]-[X_2,X_4],[X_2,X_5]-[X_3,\
X_4])} \leqno(4.1)
$$
where $ k<X_1,X_2,X_3,X_4,X_5>$ is the free associative algebra in the five variables $X_i$ and
$[X_i,X_j]=X_iX_j-X_jX_i$ is the commutator.
The corresponding Hilbert series is:
$$
{1\over R_{20}^!(z)}=\prod_{n=1}^\infty (1-z^{2n-1})^5(1-z^{2n})^3
$$
The proof of this last statement is by an adaption of the proof given in [L\"o-Ro 2] when the ring $R$
is commutative (the proof is even easier in the skew-commutative case).
\bigskip
Here is Case 12:
$$
R_{12}={E(x_1,x_2,x_3,x_4,x_5)\over(x_1x_2,x_1x_3+x_2x_4+x_3x_5,x_4x_5)}
$$
and the Hilbert series is still $R_{12}(z)=1+5z+7z^2$,
But the corresponding Hilbert series for the Koszul dual $R_{12}^!$
is given by
$$
{1\over R_{12}^!(z)}=(1-2z)^2\prod_{n=1}^\infty (1-z^n)
$$
This is proved in the same way as it was proved for the corresponding $R$ in the commutative case
 [L\"o-Ro 1].
Finally here is Case 15:
$$
R_{15}={E(x_1,x_2,x_3,x_4,x_5)\over(x_1x_4+x_2x_3,x_1x_5,x_3x_4+x_2x_5)}
$$
which has $R_{15}(z) = 1+5z+7z^2$, but
$$
{1\over R_{15}^!(z)}=(1-2z)\prod_{n=1}^\infty (1-z^{2n-1})^3(1-z^{2n})^2
$$
which is proved in a similar way.
\bigskip
Now if $R$ is any of the  3 rings above,
then then Gorenstein ring  $G=R\propto I(k)$ has
Hilbert series $G(z)= 1+12z+12z^2+z^3$
and by the Gulliksen formula
$P_G(x,y)={P_R(x,y) \over 1-xyExt^*_R(k,R)(x,y)}$,
the B{\o}gvad formula $Ext^*_R(k,R)(x,y)/P_R(x,y) = x^2y^2R(-{1\over xy})$ (cf. section 5 below)
and the L\"ofwall formula ${1\over P_R(x,y)} = (1+1/x)/R^!(xy)-R(-xy)/x$
we finally find that the transcendental properties of $P_G(x,y)$ are
rationally related to those of $R^!(z)$. For the validity of the B{\o}gvad formula,
cf the last section.
\bigskip
We now present the 3 other cases of $R$ with 4 quadratic relations (cases 21,22 and 33) where the
Hilbert series is $R(z)=1+5z+6z^2$ and which lead to possibly strange Gorenstein rings and 3-manifolds $X$ with the
dimension of $H^1(X,{\bf Q})$ equal to 11. They are extremely easy to describe:
$$
R_{21} = R_{20}/(x_4x_5), \quad R_{22} = R_{20}/(x_3x_5) \quad {\rm and}\quad R_{33}=R_{15}/(x_4x_5)
$$
but the corresponding series $R^!(z)$ are unknown. But, using the Backelin et al programme BERGMAN [Ba]
we have calculated these series up to degrees 25, 14, and 32 respectively.
In the last case we found using a work-station with 48 GB of internal memory that
$$
{1\over(1-z)^2R^!_{33}(z)} = 1-3z-z^2+z^3+2z^4+3z^5+z^6+z^7-z^8-z^9-2z^{10}-z^{11}-3z^{12}-z^{13}-z^{14}
$$
$$
-z^{15}+z^{17}+z^{18}+2z^{19}+z^{20}+z^{21}+3z^{22}+z^{23}+z^{25}+z^{26}-z^{29}-z^{30}-z^{31}-z^{32} \ldots
$$
But to go from degree 31 to degree 32 we needed more than one week of calculations.
But we still think that the series is transcendental here.
\medskip
REMARK 4.1.--- All these results are in the case the characteristic of the base field is 0.
In Case 20 we have different $R^!_{20}(z)$ for all characteristics and the same remarks seems to be applicable
to the cases 21,22 and 33. What this gives for the corresponding fundamental groups of the corresponding
3-manifolds $X$ has not been studied.
\medskip
REMARK 4.2.--- If the Yoneda $Ext$-algebra $Ext^*_R(k,k)$ is not finitely generated as an algebra
 then $Ext^*_{R \propto M}(k,k)$ is not so either, since the algebra $Ext^*_R(k,k)$ is a quotient
of  $Ext^*_{R \propto M}(k,k)$.
One can use this for the Gorenstein ring $G_{33}=R_{33}\propto I(k)$ which has Hilbert series
$G_{33}(z) = 1+11z+11z^2+z^3$. 
Now $Ext^*_{R_{33}}(k,k)$ needs an infinite number of generators if and only if 
the $Tor_{3,j}^{R^!_{33}}(k,k)$ is non-zero for an infinite number of $j:s$ (cf. e.g. Theorem 3.1, (a)  in [Ro 2]).  
Indeed, by computer calculations using the ANICK command in the programme BERGMAN one obtains that
the dimensions of $Tor_{3,j}^{R^!_{33}}(k,k)$ are 1 for $j=4$ and then $0,3,0,2,1,2$ for $j = 5,6,7,8,9,10$
and again $0,3,0,2,1,2$ for $j = 11,12,13,14,15,16$ etc. 

\mysec{5. A formula of B{\o}gvad and $m^2$-selfinjective rings.}
\bigskip
Recall that B{\o}gvad proved the following in [B{\o}]:
Let $(R,m)$ be a local commutative ring with $m^3 =0$. Assume that
$R$ has some special properties ($soc(R) = m^2$ and $R$ being a ``beast'').
 Then we have the formula for the Bass series $Bass_R(Z)$, i.e. the generating series
of $Ext^*_R(k,R)$ and the Poincar\'e-Betti series $P_R(Z)$, i.e. the generating
series of $Ext^*_R(k,k)$:
$$
Bass_R(Z)/P_R(Z)= Z^2R(-{1\over Z})\, {\rm where}\, R(Z)=1+|m/m^2|z+|m^2|Z^2\, {\rm Hilbert\, series\, of}\, R
$$
In [B{\o}] this formula is broken up into two assertions, of which the first one is often valid
(the proposition ``$E(R/m^2)$'') and the other is more special.
Lescot has observed in [Les 2], [Les 4] that the proof of this in [B{\o}] boils 
down to prove that CONDITION 5.ii and CONDITION 5.iii below are valid  under some conditions:
\bigskip
CONDITION 5.i --- The natural map $Ext^*_R(k,m)\rightarrow Ext^*_R(k,R)$ is an epimorphism.
\medskip
CONDITION 5.ii --- The natural map $Ext^*_R(k,m/m^2)\rightarrow Ext^*_R(k,R/m^2)$ is an epimorphism.
\medskip
CONDITION 5.iii --- The natural map $Ext^*_R(k,m^2)\rightarrow Ext^*_R(k,R)$ is an epimorphism (for $*=0$ this means
that the socle of $R$ is $m^2$)
\bigskip
The CONDITION 5.i is true if $(R,m)$ is nonregular. The condition CONDITION 5.ii
is also often true (cf. condition ``$E(R/m^2)$'' in [B{\o}] and and Proposition 1.10 in [Les 2] for 
$I=m^2$, as well as the assertion that $G \rightarrow G/soc(G)$ is a Golod
map for a Gorenstein ring $G$ [Av-Le].
We have not yet proved all skew-commutative versions of the preceding
results, but computer computations indicate that they are true up
to ``high degrees'' and probably in all degrees.
\bigskip
We now present some results that should give the skew-versions of some of the
preceding results. We hope to return to these problems rather soon.
\medskip
Let $R$ be any ring (with unit). Recall that in order to test that a left $R$-module $M$ is injective it is
sufficient to test that for any left $R$-ideal $J$,  any $R$-module map
$\phi:  J \rightarrow M$ can be extended to a map $R \rightarrow M$. Since the last map is given as 
$r \rightarrow r.m_{\phi}$ where $m_{\phi}$ is a suitable element of $M$,
 this means that $\phi(j)-jm_{\phi}=0$ for any $j\in J$. In particular $R$ is self-injective
to the left if and only if for any $R$-module map $\phi: J \rightarrow R$ there is an element $r_{\phi} \in R$
such that $\phi(j)- j.r_{\phi}= 0$ for all $j\in J$. This explains the condition b) in the Lemma that follows:
\bigskip
LEMMA 5.1.--- Let $(R,m)$ be a local (skew)commutative local ring where $m^3=0$ and $J$ an ideal in $R$.
 The following two conditions are equivalent:

a) The natural map
$$
Ext^1_R(R/J,m^2) \longrightarrow Ext^1_R(R/J,R)
$$
is surjective.

b) For any $R$-module map $\phi: J \rightarrow R$ there is an element $r_{\phi} \in R$ such that
$\phi(j)- j.r_{\phi} \in m^2$ for all $j\in J$.
\bigskip
PROOF: The short exact sequence $0\rightarrow J \rightarrow R \rightarrow R/J \rightarrow 0$
and the natural map $m^2 \rightarrow R$ give rize to a commutative diagram with exact rows:
$$
 0 \rightarrow Hom_R(R/J,m^2)\rightarrow Hom_R(R,m^2)\rightarrow Hom_R(J,m^2)\buildrel{\delta'}\over\rightarrow Ext^1_R(R/J,m^2)\rightarrow  0
$$
$$
    \downarrow \quad \quad \quad \quad  \quad \quad \quad \quad \downarrow \quad \quad \quad \quad \quad \quad \quad \quad \downarrow \kappa \quad \quad \quad \quad \quad \quad \quad \quad \quad \downarrow j
$$
$$
0 \rightarrow Hom_R(R/J,R)\rightarrow Hom_R(R,R)\buildrel{i}\over\rightarrow Hom_R(J,R)\buildrel{\delta}\over\rightarrow Ext^1_R(R/J,R)\rightarrow  0
$$
Let us first prove that  $ a) \Rightarrow b)$. Thus assume that $j$ is onto. We start with
a $\phi \in Hom_R(J,R)$ and put $\alpha = \delta(\phi)$.
We can assume that $\alpha = j(\xi)$. But $\xi = \delta'(\tilde\xi)$.
Thus $\delta(\kappa(\tilde\xi)-\phi)=0$, so that $\kappa(\tilde\xi)-\phi$ comes by $i$ from a map $R \rightarrow R$
of the form $ r \rightarrow r.r_{\phi} $ . Thus to any map $\phi : J \rightarrow R$ there is a $r_{\phi}$ in $R$
such that for all $j\in J$, $\eta(j)-j.r_{phi} \in m^2$, i.e. we have proved b). The converse follows from the ``reverse''
reasoning.
\bigskip
REMARK 5.1.--- If $L$ is any $R$-module and $P(L) \rightarrow L$ is the projective envelope of $L$ we have an exact
sequence
$$
0 \rightarrow S(L) \rightarrow P(L) \rightarrow L \rightarrow 0
$$ 
where the first syzygy $S(L)$ of $L$ is included in $m.P(L)$, we can redo the same reaoning for $P(L)/S(L)$ as we did for
$R/J$ in LEMMA 5.1. The result is that $Ext^1_R(L,m^2) \rightarrow Ext^1_R(L,R)$ is onto
if and only if for any map $\phi : S(L) \rightarrow R$ there is a map $j_{\phi}: P(L) \rightarrow R$ such that
$\phi(s)-j_{\phi}(s) \in m^2$ for all $s\in S(L) \subset P(L)$ Note that $P(L)$ is free so 
that $j_{\phi}$ is given by a matrix
of elements in $R$. 
\bigskip
REMARK 5.2.--- In the selfinjective case ($0$-selfinjective ) it is of course sufficient to require b) for the maximal ideal $J=m$.
In the ``$m^2$-selfinjective'' case it is not clear (for us) what the right definition should be. We hope to return to this later. Therefore the definition below is
maybe too strong.
\bigskip
DEFINITION 5.1.--- We say that $R$ is $m^2$-selfinjective if the conditions of REMARK 5.2 are valid
for $L=k$ and all syzygies of $k$.
\bigskip
Using this definition we can formulate the following consequence of LEMMA 5.1 and and REMARK 5.1:
\bigskip
COROLLARY 5.1.--- The following conditions are equivalent:

$\alpha$) $Ext^*_R(k,m^2) \rightarrow Ext^*_R(k,R)$ is onto
\medskip
$\beta$)  $R$ is $m^2$-selfinjective.
\bigskip
REMARK 5.3. --- In the commutative case, it seems that many rings $(R,m)$ with $m^3=0$
of the form $R=k[x_1,x_2,\ldots,x_n]/(f_1,f_2,\ldots,f_t)$ where
the $f_i$ are homogeneous quadratic forms are $m^2$-selfinjective.
For the case when $n=4$ cf. REMARK 5.4 below. 
\bigskip
OBSERVATION 5.1.--- Let $(R,m)$ be a local ring with $m^3=0$ residue field $k=R/m$ and $soc(R)=m^2$.

The following two conditions are equivalent:

a) $R$ is $m^2$-selfinjective and condition CONDITION 5.ii above is true.

b) The Bass series of $R$, i.e.
$Bass_R(Z)=\sum_{i\geq 0} |Ext^i(k,R)| Z^i$ is related to the Poincar\'e-Betti series 
$P_R(Z) = \sum_{i\geq 0} |Ext^i_R(k,k)| Z^i$ by the ``B{\o}gvad formula''
$Bass_R(Z) = Z^2R(-1/z)P_R(Z)$,
where $R(Z) = 1+|m/m^2|Z+|m^2|Z^2$ is the Hilbert series of $R$.

PROOF:- Consider the long exact sequence obtained when we apply $Ext^*_R(k,-)$
to the short exact sequence $0 \rightarrow m^2 \rightarrow R \rightarrow R/m^2 \rightarrow 0$:
$$
0 \rightarrow Hom_R(k,m^2)\rightarrow Hom_R(k,R)\rightarrow Hom_R(k,R/m^2)\rightarrow Ext^1_R(k,m^2) \rightarrow Ext^1_R(k,R)\rightarrow
$$
$$
 \rightarrow Ext^1_R(k,R/m^2)\rightarrow Ext^2_R(k,m^2)\rightarrow Ext^2_R(k,R)\rightarrow Ext^2_R(k,R/m^2)\rightarrow 
Ext^3_R(k,m^2) \rightarrow 
$$
Since $soc(R)=m^2$ the natural monomorphism $Hom_R(k,m^2)\rightarrow Hom_R(k,R)$ is indeed an isomorphism.
Now according to a) all the maps $Ext^i_R(k,m^2) \rightarrow Ext^i_R(k,R)$ are also epimorphisms for $i \geq 1$.
It follows that we have an exact sequence:
$$
 0 \rightarrow s^{-1}Ext^*_R(k,R/m^2)\rightarrow Ext^*_R(k,m^2)\rightarrow Ext^*_R(k,R)\rightarrow 0
$$
so that $Bass_R(Z)-|m^2|P_R(Z)+Z.Ext^*_R(k,R/m^2)(Z)=0$.
But it is now easy to apply CONDITION 5.ii and this gives the result that $a)\Rightarrow b)$.
The converse is easy.
\bigskip
REMARK 5.4.--- In the commutative case the only rings of embedding dimension 4 which are quotients of $k[x,y,z,u]$ by an ideal
$I$ generated by homogeneous quadratic forms are and having $m^3=0$ are according to [Ro 5]
given by (the numbers of the ideals comes from [Ro 5]:
\medskip
$I_{29} = (x^2+xy,y^2+xu,z^2+xu,zu+u^2,yz)$
\medskip
$I_{54} = (x^2,xy,y^2,z^2,yu+zu,u^2)$
\medskip
$I_{55} = (x^2+xy,xz+yu,xu,y^2,z^2,zu+u^2)$
\medskip
$I_{56} = (x^2+xz+u^2,xy,xu,x^2-y^2,z^2,zu)$
\medskip
$I_{57} = (x^2+yz+u^2,xu,x^2+xy,xz+yu,zu+u^2,y^2+z^2)$
\medskip
$I_{71} = (x^2,y^2,z^2,u^2,xy,zu,yz+xu) $
\medskip
$I_{78} = (x^2,xy,y^2,z^2,zu,u^2,xz+yu,yz-xu)$
\medskip
$I_{81} = (x^2,y^2,z^2,u^2,xy,xz,yz-xu,yu,zu)$
\medskip
The Hilbert series $R(z)$ of the different cases $R=k[x,y,z,u]/I$ are
$$
1+4z+5z^2\quad {\rm (case\,\, 29)},\quad 1+4z+4z^2\quad {\rm (cases\,\, 54,55,56,57)},
$$
$$
1+4z+3z^2\quad {\rm (case\,\, 71)},\quad 1+4z+2z^2\quad {\rm (case\,\, 78)\quad and}\quad 1+4z+z^2\quad {\rm (case\,\, 81)} 
$$
But the Hilbert series $R^!(z)$ of the Koszul duals $R^!$ are respectively (all different):
$$
{(1+z)^4\over(1-z^2)^5},\quad {1\over 1-4z+4z^2},\quad {(1-z+z^2)^2\over (1-z)^3(1-3z+3z^2-3z^3)},\quad {1-z+z^2 \over (1-z)^2(1-3z+2z^2-z^3)} 
$$
$$
{1 \over (1-z)^2(1-2z-z^2)},\quad {1 \over 1-4z+3z^2},\quad  {1 \over 1-4z+2z^2}\quad{\rm and}\quad  {1 \over 1-4z+z^2}
$$
and in all cases the ``L\"ofwall formula''
$$
{1\over P_R(x,y)}= (1+{1/x})/R^!(xy)-R(-xy)/x
$$
holds true.
In all these cases {\it except} $I_{78}$ the B{\o}gvad formula also holds true.

\bigskip
REMARK 5.5.--- The B{\o}gvad formula is indeed a two-variable formula:
$$
Bass_R(x,y)/P_R(x,y)= x^2y^2R(-{1\over xy})
$$
 which shows that the non-diagonal elements occur for the two-variable
Bass series in ``the same way'' as they occur in the two-variable Poincar\'e-Betti series.
In the case of $I_{78}$ this is not true. Indeed $I_{78}$ is a Koszul ideal but the corresponding
Bass series has non-diagonal elements. More precisely we have in that case:
$$
Bass_R(x,y)/P_R(x,y) = x^2y^2R(-{1\over xy})+x^2y^2+xy^2
$$
REMARK 5.6.--- There are of course similar results to those of REMARK 5.4 and 5.5 in the skew-commutative case.
What corresponds to the ``bad'' case $I_{78}$ in the case of four commuting variables is the case
of four skew-commuting variables: $E(x,y,z,u)\over (xy,xz-yu,yz-xu,zu)$.
\bigskip
REMARK 5.7.--- In general it is not true that $Bass_R(x,y)/P_R(x,y)$ is a rational function.
This was first noted by Lescot in his thesis [Les 4]. Take e.g. $S$ a ring with transcendental Bass series.
Form $T = S \propto I(S/m)$ and let $R$ be $T/(socle T)$ then $R$ is a so-called Teter ring [Tet], where
the maximal ideal is isomorphic to its Matlis dual (this even characterizes Teter rings [Hu]), and from
this one sees that $Bass_R(x,y)/P_R(x,y)$ is transcendental [Les 2,Corollaire 1.9 ]. 

\bigskip
REFERENCES
\bigskip
[An] Anick, David J.; Gulliksen, Tor H.
Rational dependence among Hilbert and Poincar\'e series. 
J. Pure Appl. Algebra 38 (1985), no. 2-3, 135--157.

[Av-Le] Avramov, Luchezar L.; Levin, Gerson L. 
Factoring out the socle of a Gorenstein ring. 
J. Algebra 55 (1978), no. 1, 74--83.

[Ba] Backelin, J\"orgen et al, BERGMAN, a programme for non-commutative Gr\"obner basis
calculations available at
{\tt http://servus.math.su.se/bergman/}

[Be] Benson, Dave
An algebraic model for chains on $\Omega BG^\wedge_p$. 
Trans. Amer. Math. Soc. 361 (2009), no. 4, 2225--2242.

[B{\o}] B{\o}gvad, Rickard
Gorenstein rings with transcendental Poincar\'e	-series. 
Math. Scand. 53 (1983), no. 1, 5--15.

[Co-Su] Cohen, Daniel C.; Suciu, Alexander I. The boundary manifold of a complex line arrangement. Groups, homotopy and configuration spaces, 105--146, Geom. Topol. Monogr., 13, Geom. Topol. Publ., Coventry, 2008.

[Elk-Sri] El Khoury, Sabine; Srinivasan, Hema
A class of Gorenstein artin algebras of embedding dimension four. Communication in Algebra, 37 (2009), 3259-3257.

[Fo-Gr-Rei] Fossum, Robert M.; Griffith, Phillip A.; Reiten, Idun
Trivial extensions of abelian categories. 
Homological algebra of trivial extensions of abelian categories with applications to ring theory. Lecture Notes in Mathematics, Vol. 456. Springer-Verlag, Berlin-New York, 1975. xi+122 pp.

[Fr\"o] Fr\"oberg, Ralph
Determination of a class of Poincar\'e	 series. 
Math. Scand. 37 (1975), no. 1, 29--39.

[Gu] Gulliksen, Tor H.
Massey operations and the Poincar\'e series of certain local rings. 
J. Algebra 22 (1972), 223--232.

[Gur] Gurevich, G. B.
Foundations of the theory of algebraic invariants. 
Translated by J. R. M. Radok and A. J. M. Spencer P. Noordhoff Ltd., Groningen 1964 viii+429 pp.

[Hu] Huneke, Craig ; Vraciu, Adela
 Rings which are almost Gorenstein
   
{\tt arXiv:math/0403306}

[I] Ivanov, A. F.
Homological characterization of a class of local rings. (Russian) 
Mat. Sb. (N.S.) 110(152) (1979), no. 3, 454--458, 472.

[Les 1] Lescot, Jack, S\'eries de Bass des modules de syzygie. (French) [Bass series of syzygy modules] Algebra, algebraic topology and their interactions (Stockholm, 1983), 277--290, 
Lecture Notes in Math., 1183, Springer, Berlin, 1986.

[Les 2] Lescot, Jack,
La s\'erie de Bass d'un produit fibr\'e	 d'anneaux locaux. (French) [The Bass series of a fiber product of local rings] Paul Dubreil and Marie-Paule Malliavin algebra seminar, 35th year (Paris, 1982), 218--239, 
Lecture Notes in Math., 1029, Springer, Berlin, 1983.

[Les 3] Lescot, Jack
Asymptotic properties of Betti numbers of modules over certain rings. 
J. Pure Appl. Algebra 38 (1985), no. 2-3, 287--298.

[Les 4]  Lescot, Jack
Contribution \`a l'\'etude des s\'eries de Bass, Th\`ese, Universit\'e de Caen, 1985.

[Lev] Levin, Gerson
Modules and Golod homomorphisms. 
J. Pure Appl. Algebra 38 (1985), no. 2-3, 299--304.

[L\"of] L\"ofwall, Clas, On the subalgebra generated by the one-dimensional elements in the Yoneda Ext-algebra. Algebra, algebraic topology and their interactions (Stockholm, 1983), 291--338, 
Lecture Notes in Math., 1183, Springer, Berlin, 1986.

[L\"of-Ro 1] L\"ofwall, Clas; Roos,Jan-Erik,
Cohomologie des alg\`ebres de Lie gradu\'ees et s\'eries de Poincar\'e-Betti non rationnelles.
 (French. English summary) 
C. R. Acad. Sci. Paris S\'er. A-B 290 (1980), no. 16, A733--A736.

[L\"of-Ro 2] L\"ofwall,Clas; Roos,Jan-Erik,
A nonnilpotent $1$-$2$-presented graded Hopf algebra whose Hilbert series converges in the unit circle. 
Adv. Math. 130 (1997), no. 2, 161--200.

[Mat] Matlis, Eben
Injective modules over Noetherian rings. 
Pacific J. Math. 8 1958 511--528.

[Mil] Milnor, J.
A unique decomposition theorem for $3$-manifolds. 
Amer. J. Math. 84 1962 1--7.

[Pos] Positselski., L. E.
The correspondence between Hilbert series of quadratically dual algebras does not imply their having the Koszul property. (Russian) Funktsional. Anal. i Prilozhen. 29 (1995), no. 3, 83--87; translation in 
Funct. Anal. Appl. 29 (1995), no. 3, 213--217 (1996)

[Ro 1] Roos, Jan-Erik
Relations between Poincar\'e-Betti series of loop spaces and of local rings. S\'eminaire d'Alg\`ebre Paul Dubreil 
31\`eme ann\'ee (Paris, 1977--1978), pp. 285--322, 
Lecture Notes in Math., 740, Springer, Berlin, 1979.

[Ro 2] Roos, Jan-Erik,
The homotopy Lie algebra of a complex hyperplane arrangement is not necessarily finitely presented. (English summary) 
Experiment. Math. 17 (2008), no. 2, 129--143.

[Ro 3] Roos, Jan-Erik, 
Homological properties of quotients of exterior algebras, In preparation,
Abstract available at {\it Abstracts Amer. Math. Soc.} 21 (2000), 50-51.

[Ro 4] Roos, Jan-Erik,
On the characterisation of Koszul algebras. Four counterexamples. (English, French summary) 
C. R. Acad. Sci. Paris S\'er. I Math. 321 (1995), no. 1, 15--20.

[Ro 5] Roos, Jan-Erik,
A description of the Homological Behaviour of Families of Quadratic Fourms in Four variables,
Syzygies and Geometry: Boston 1995 (A.Iarrobino, \break A.Martsinkovsky and J.Weyman,editors),
Northeastern Univ. 1995, 86-95

[Sik] Sikora, Adam S.
Cut numbers of 3-manifolds. (English summary) 
Trans. Amer. Math. Soc. 357 (2005), no. 5, 2007--2020 (electronic).

[Sul] Sullivan, Dennis
On the intersection ring of compact three manifolds. 
Topology 14 (1975), no. 3, 275--277.

[Tet] Teter, William
Rings which are a factor of a Gorenstein ring by its socle. 
Invent. Math. 23 (1974), 153--162.

[Vin-El] Vinberg,B.;Elashvili, A. G.

A classification of the three-vectors of 
nine-dimensional space. (Russian)

Trudy Sem. Vektor. Tenzor. Anal. 18 (1978), 197--233.

(English translation in: Selecta Math. Soviet. 7 (1988), no. 1. 63-98.)
\par}
\end